\documentclass[11pt]{article}
\usepackage{a4}
\usepackage{amsmath}
\usepackage{amssymb}

\font\tenmsy=msbm10 scaled 1200 \font\sevenmsy=msbm7 scaled 1200
\font\fivemsy=msbm5 scaled 1200
\newfam\msyfam
\textfont\msyfam=\tenmsy \scriptfont\msyfam=\sevenmsy
\scriptscriptfont\msyfam=\fivemsy

\newcommand{\R}{{\mathbb R}}

\newcommand{\Z}{{\mathbb Z}}
\newcommand{\N}{{\mathbb N}}

\newcommand{\Half}{{\mathbb H}}

\newcommand{\cH}{{\cal H}}
\newcommand{\hf}{{\cal H}^{f}}
\newcommand{\hs}{{\cal H}^{s}}

\newcommand{\p}{\psi}

\newcommand{\ep}{ \epsilon }

\newcommand{\Ja}{Jarn\'{\i}k }

\newcommand{\QED}{{\hspace*{\fill}\raisebox{-1ex}{$\spadesuit$}}}

\newtheorem{thschmidt}{Theorem BDV1 }
\newtheorem{thdet}{Theorem BDV2 }
\newtheorem{thkh}{Khintchine's Theorem (1924) }
\newtheorem{thjar}{Jarn\'{\i}k's Theorem (1931) }

\newtheorem{thbs}{Schmidt's Theorem (1964).}

\newtheorem{problem}{Problem}
\newtheorem{lemma}{Lemma}
\newtheorem{theorem}{Theorem}
\newtheorem{corollary}{Corollary}
\newtheorem{proposition}{Proposition}
\newtheorem{adeq}{Theorem (Khintchine--Jarn\'{\i}k)}

\newcommand{\cR}{{\cal R}}
\newcommand{\ka}{\kappa}
\newcommand{\ra}{R_{\alpha}}
\newcommand{\De}{\Delta}
\newcommand{\ma}{\beta_\alpha}
\newcommand{\dist}{\mathrm{dist}\,}
\renewcommand{\r}{\rho}
\newcommand{\La}{\Lambda}

\newcommand{\be}{\begin{eqnarray*}}
\newcommand{\ee}{\end{eqnarray*}}

\newcommand{\pen}{{\rm pen}}

\begin{document}

%\title{\huge Ubiquity, a general logarithm law for geodesics and Diophantine approximation
%\\
%on fractal sets.
%}

\title{\huge Ubiquity and a general logarithm law  \\ for geodesics. }

\author{Victor Beresnevich\footnote{EPSRC Advanced Research Fellow, EP/C54076X/1}
\\ {\small\sc York} \and Sanju Velani%\footnote{Royal Society University Research
%Fellow}
\\ {\small\sc York}}

%\author { Sanju Velani\footnote{Royal Society University Research
%Fellow} \\ {\small Y{\scriptsize ORK}} }

\date{ In memory of Bill Parry \\ ~ \\
} \maketitle

%{\small
\begin{abstract}

There are two fundamental results in the classical theory of
metric Diophantine approximation:   Khintchine's theorem  and
Jarn\'{\i}k's theorem. The former relates the size of the set of
well approximable numbers, expressed in terms of Lebesgue measure,
to the behavior of a certain volume sum. The latter is a Hausdorff
measure version of the former. We start by discussing these
theorems and  show  that they are both in fact a simple
consequence of the notion of `local ubiquity'. The local ubiquity
framework introduced here is a much simplified and  more
transparent version of that in \cite{memoirs}. Furthermore, it
leads to a  single local ubiquity theorem that unifies the
Lebesgue and Hausdorff theories. As an application of our
framework we consider the theory of metric Diophantine
approximation on limit sets of Kleinian groups. In particular, we
obtain a general Hausdorff measure version of Sullivan's
logarithm law for geodesics -- an aspect overlooked in \cite{memoirs}.  %Secondly, we consider the problem of
%Diophantine approximation on fractal sets -- in particular, on the
%standard middle third Cantor set $K$ and $M_o$ sets. Our results
%for $K$ are shown to answer a problem of K. Mahler. As with all
%`good' problems its solution opens up a can of worms.

\vspace{5mm}

%\noindent {\bf 2000} {\em  Mathematics Subject Classification:}
%11J83, 11J13, 11K60, 28A78, 28A80
\end{abstract}
%}

%\newpage

\section{Introduction}

\subsection{Background: the classical theory\label{beg}}
To set the scene, we follow the opening discussion of
\cite{memoirs} and introduce a basic $\limsup$ set whose study has
played a central role in the development of the classical theory
of metric Diophantine approximation. Given a real, positive
decreasing function $\psi : \R^+ \to \R^+$, let $$ W (\p) := \{x
\in [0,1] :|x - p/q| < \p(q) \; {\rm for\ i.m.\ rationals\ }  p/q
\   (q>0) \}, $$ where `i.m.' means `infinitely many'. This is the
classical set of $\p$--well approximable numbers in the theory of
Diophantine approximation. The fact that we have restricted our
attention to the unit interval rather than the real line is purely
for convenience. It is natural to refer
 to the  function $\p$ as the {\em approximating function}.
It governs the `rate' at which points in the unit interval must be
approximated by rationals  in order to lie in $W(\p)$. It is not
difficult to see that $ W (\p) $ is a $\limsup$ set. For $n \in
\N$, let $$
 W(\p,n) :=  \!\!\!\!\!  \bigcup_{ { k^{n-1} < q  \leq k^n}} \bigcup_{ { 0 \leq  p  \leq q}}
\!\! B(p/q ,\p(q))  \cap [0,1]
$$ where  $k>1$ is
fixed and $B(c,r)$ is the open interval centred at $c$ of radius
$r$.  The set $ W(\psi)$ consists precisely  of points in the unit
interval that lie in infinitely many $ W(\p,n)$; that is
$$ W(\psi) = \limsup_{n \to \infty} W(\p,n) :=
\bigcap_{m=1}^{\infty} \bigcup_{n=m}^{\infty} W(\p,n) \ . $$

Investigating the measure theoretic properties of the set $W(\p)$
underpins the classical theory of metric Diophantine
approximation. We begin by considering the `size' of $W(\p)$
expressed in terms of the ambient measure $m$; i.e.
one--dimensional Lebesgue measure. On exploiting the $\limsup $
nature of $W(\p)$, a  straightforward application of the
convergence part of the Borel--Cantelli lemma from probability
theory yields that $$ m(W (\p)) = 0  \ \ \ \ \ {\rm if \ } \ \ \ \
\ \sum_{n=1}^{\infty} k^{2n} \p(k^n) \ < \ \infty \ . $$ Notice
that since $\psi$ is monotonic, the convergence/divergence
property of the above sum is equivalent to that of
$\sum_{r=1}^{\infty} r \, \p(r) $.

A natural problem now arises.  Under what conditions is  $m(W
(\p)) >0 $   ?  The  following fundamental result provides a
beautiful and simple criteria for the `size' of the set $W(\p)$
expressed in terms of Lebesgue measure.

\bigskip

\begin{thkh} Let $\p$ be a real, positive decreasing function.
Then
$$ m(W(\p)) =\left\{
\begin{array}{ll}
0 & {\rm if} \;\;\; \sum_{r=1}^{\infty} \; r \,  \p(r)  <\infty\;
,\\[1ex]
 1 & {\rm if} \;\;\; \sum_{r=1}^{\infty} \; r \,  \p(r)
 =\infty \; .
\end{array}\right.$$
\end{thkh}

Thus, in the divergence case,  which constitutes the main
substance of Khintchine's theorem,  not only do we have positive
Lebesgue measure but full Lebesgue measure. To the best of our
knowledge, this turns out to be the case for all naturally
occurring limsup sets -- not just within the number theoretic
setup. Usually, there is a standard argument which allows one to
deduce full measure from positive measure -- such as the
invariance of the $\limsup$ set or some related set, under an
ergodic transformation. In any case, we shall prove a general
result which directly implies the above full measure statement.
{\em It is worth mentioning that in Khintchine's original
statement the stronger hypothesis that $r^2 \p(r)$ is decreasing
was assumed.} The fact that this additional hypothesis is
unnecessary has been known for sometime.

Returning to the convergence  case,  we cannot obtain any further
information regarding the `size' of $W(\p)$ in terms of Lebesgue
measure  --- it is always zero. Intuitively, the `size' of $W(\p)$
should decrease as the rate of approximation governed by the
function $\p$ increases. In short, we require a more delicate
notion of `size' than simply Lebesgue measure. The appropriate
notion of `size' best suited for describing the finer measure
theoretic structures of $W(\p)$ is that of generalized Hausdorff
measures. The Hausdorff $f$--measure %$\mathcal{ H}^{f}$
with respect to a dimension function $f$ is  a natural
generalization of Lebesgue measure. A {\bf dimension function} $f
: \R^+ \to \R^+ $ is an increasing, continuous  function such that
$f(r)\to 0$ as $r\to 0 \, $. The Hausdorff $f$--measure with
respect to the dimension function $f$ will be denoted throughout
by $\mathcal{ H}^{f}$ and is defined as follows. Suppose $F$ is  a
non--empty subset  of a metric space  $(\Omega,d)$. For $\rho >
0$, a countable collection $ \left\{B_{i} \right\} $ of balls in
$\Omega$ with radii $r_i  \leq \rho $ for each $i$ such that $F
\subset \bigcup_{i} B_{i} $ is called a {\em $ \rho $-cover}\/ for
$F$. For a dimension function $f$ define $ \mathcal{ H}^{f}_{\rho}
(F) \, = \, \inf \left\{ \sum_{i} f(r_i) \ : \{ B_{i} \}  {\rm \
is\ a\ } \rho {\rm -cover\ of\ } F \right\} \,
 $, where the infimum is over all $\rho$-covers.  The {\bf
Hausdorff $f$--measure} $ \mathcal{ H}^{f} (F)$ of $F$ with
respect to the dimension function $f$ is defined by   $$ \mathcal{
H}^{f} (F) := \lim_{ \rho \rightarrow 0} \mathcal{ H}^{f}_{\rho}
(F)  \; = \; \sup_{\rho > 0 } \mathcal{ H}^{f}_{\rho}  (F) \; . $$
In the case that  $f(r) = r^s$ ($s \geq 0$), the measure $ \hf $
is the usual {\em $s$--dimensional Hausdorff measure}\/ $\hs $ and
the {\bf Hausdorff dimension}  $\dim F$ of a set $F$ is defined by
$ \dim \, F \, := \, \inf \left\{ s : \mathcal{ H}^{s} (F) =0
\right\} = \sup \left\{ s : \mathcal{ H}^{s} (F) = \infty \right\}
$. In particular when $s$ is an integer $\hs$ is a constant
multiple of $s$--dimensional Lebesgue measure. For further details
see \cite{falc,jh,mat}.

%
%We referee the reader to \S\ref{hmd} for the  standard definition
%of $\mathcal{ H}^{f}$ and further comments regarding Hausdorff
%measures and dimension.

Again on exploiting the $\limsup $ nature of $W(\p)$, a
straightforward covering argument provides a simple convergence
condition under which $\mathcal{ H}^{f}(W(\p)) = 0 $. Thus, in
view of  the development of the Lebesgue theory it is  natural to
ask for conditions under which $\mathcal{ H}^{f}(W(\p))$  is
strictly positive.

The  following fundamental result provides a beautiful and simple
criteria for the `size' of the set $W(\p)$ expressed in terms of
Hausdorff measures.

\begin{thjar}  \label{main}
 Let $f$ be a dimension function such that
$r^{-1} \, f(r)\to \infty$ as $r\to 0 \, $ and  $r^{-1} \, f(r) $
is decreasing. Let $\p$ be  a real, positive  decreasing function.
Then $$ \hf\left(W(\p)\right)=\left\{\begin{array}{cl} 0 & {\rm \
if} \;\;\; \sum_{r=1}^{\infty}  \ \; r \,  f\left(\p(r)\right)
 <\infty \; ,\\[1ex]
\infty & {\rm \ if} \;\;\; \sum_{r=1}^{\infty} \  \; r \,
f\left(\p(r)\right) =\infty \; .
\end{array}\right.$$
\end{thjar}

Clearly  the above theorem  can be regarded as the Hausdorff
measure version of Khintchine's theorem.  As with the latter, the
divergence part constitutes the main substance. Notice, that the
case when $ \hf $ is comparable to one--dimensional Lebesgue
measure $m$ (i.e. $ f(r)= r$) is excluded by the condition $r^{-1}
\, f(r)\to \infty$ as $r\to 0 \, $. Analogous to Khintchine's
original statement, {\em in Jarn\'{\i}k's original statement the
additional hypotheses that $r^{2}\p(r)$ is decreasing, $r^{2}\p(r)
\to 0 $ as $ r \to \infty $ and that $r^2f(\p(r)) $ is decreasing
were assumed.} Thus, even in the simple case when $f(r) = r^s $ $
(s\geq 0) $  and the approximating function is given by $\psi(r) =
r^{-\tau} \log r $ $(\tau > 2)$,  Jarn\'{\i}k's original statement
gives no information regarding the $s$--dimensional Hausdorff
measure  of $W(\p)$ at the critical exponent $s=2/\tau$ -- see
below. That this is the case is due to the fact that  $r^2f(\p(r))
$ is not decreasing. However, as we shall see these additional
hypotheses are unnecessary. Furthermore, with the theorems of
Khintchine and Jarn\'{\i}k as stated above it is possible to
combine them to obtain a \textbf{single unifying statement} (see
\S \ref{basiceg}) that provides a complete measure theoretic
description of $W(\p)$.

Returning to Jarn\'{\i}k's theorem,  note that in the case when $
\hf $ is the standard $s$--dimensional Hausdorff measure $\hs$
(i.e. $ f(r)= r^s $), it follows from the definition of Hausdorff
dimension
%(see \S\ref{hmd})
that $$ \dim W(\p) \, =  \, \inf \{ s
: \mbox{$\sum_{r=1}^\infty $} \; r \, \p(r)^s  < \infty \} \; . $$

Previously, Jarn\'{\i}k (1929) and independently Besicovitch
(1934) had determined the Hausdorff dimension of the set
$W(r\mapsto r^{-\tau})$, usually denoted by $W(\tau)$, of
$\tau$--well approximable numbers. They proved that for $\tau >
2$, $\dim W(\tau) = 2/\tau $. Thus, as the `rate' of approximation
increases (i.e. as $\tau$ increases) the `size' of the set
$W(\tau)$ expressed in terms of Hausdorff dimension decreases. As
discussed earlier, this is in precise keeping with one's
intuition. Obviously, the dimension result implies that $$
\mathcal{ H}^s \left(W(\tau)\right)=\left\{\begin{array}{ll} 0 & \
\ \ \
{\rm if} \;\;\; s \; > \; 2/\tau   \\
%& \\
\infty & \ \ \ \ {\rm if} \;\;\; s \; < \; 2/\tau
\end{array}\right. ,$$
but gives no information regarding the $s$--dimensional Hausdorff
measure of $W(\tau)$ at the critical value $s=\dim W(\tau)$.
Clearly, Jarn\'{\i}k's zero--infinity law implies the dimension
result and that for $\tau > 2$ $$ \cH^{2/\tau}(W(\tau)) \, = \,
\infty \ . $$ Furthermore, the `zero--infinity' law allows us to
discriminate between sets with the same dimension and even the
same $s$--dimensional Hausdorff measure. For example, with $\tau
\geq 2 $ and   $ 0 < \ep_1 < \ep_2 $ consider the approximating
functions $$ \p_{\ep_i}(r) \, :=  \, r^{-\tau} \, \left( \log \, r
\right)^{ -  \frac{\tau}{2} \left(1 + \ep_i \right) } \hspace{1cm}
(i=1,2)  \; \; . $$ It is easily verified that for any $\ep_i >
0$, $$ m(W(\p_{\ep_i})) = 0 \, , \quad \dim W(\p_{\ep_i}) = 2/\tau
\quad {\rm and}\quad \cH^{2/\tau} (W(\p_{\ep_i})) =0 \ . $$
However, consider the dimension  function $f$ given by $ f(r) =
r^{2/\tau}( \log r^{-1/\tau})^{\ep_1}$.  Then $ \sum_{r=1}^\infty
r \, f\left(\p_{\ep_i} (r)\right) \asymp \sum_{r=1}^\infty \; (r
\, (\log r)^{1+ \ep_i - \ep_{1} } )^{-1} $, where as usual the
symbol $\asymp $ denotes comparability (the quotient of the
associated quantities is bounded from above and below by positive,
finite constants). Hence, Jarn\'{\i}k's zero--infinity law implies
that
\[
\hf\left(W(\p_{ \ep_1 } )\right) \ = \ \infty \hspace{12mm}  {\rm
whilst \ }  \hspace{12mm}    \hf\left(W(\p_{ \ep_2 } )\right) \ =
\ 0 \; .
\]
Thus the Hausdorff measure $ \hf $ does make a distinction between
the `sizes' of the sets  under consideration; unlike
$s$--dimensional Hausdorff measure.

Within this classical setup, it is apparent that Khintchine's
theorem together with Jarn\'{\i}k's zero--infinity law provide a
complete measure theoretic description of $W(\p)$  -- see \S
\ref{basiceg} for a single unifying statement. In short, our
central  aim is to establish analogues of the divergence parts of
these classical results within a general framework. Recall, that
the divergence parts constitute the main substance of the
classical statements.

\subsection{The general setup and fundamental problems \label{gensetup}}
The setup described below is a much simplified version of that
considered in \cite{memoirs}. In particular, we make no attempt to
incorporate the linear forms theory of metric Diophantine
approximation. However this does have the advantage of making the
exposition more transparent and also leads to cleaner statements
which are more than adequate for the application we have in mind.

Let $(\Omega,d)$ be a compact metric space equipped with a
non-atomic, probability measure $m$. Let $\cR=\{\ra \subset \Omega
: \alpha \in J \}$  be a family of points $\ra$ of $\Omega$
indexed by an infinite, countable set $J$. The points $\ra$ will
be referred to as {\em resonant points} for reasons which will
become apparent later. Next, let $\beta: J \to \R^+ : \alpha
\mapsto \ma $ be a positive function on $J$. Thus, the function
$\beta$ attaches a `weight' $\ma$  to the resonant point $\ra$. To
avoid pathological situations within our framework, we shall
assume that the number of $\alpha$ in $J$ with $\ma$ bounded above
is always finite.

 Given a
decreasing function $\psi : \R^+ \to \R^+ $ let
$$\La(\psi)=\{x\in\Omega:x\in B(\ra,\p(\ma))\ \mbox{for\
infinitely\ many\ }\alpha\in J\} \ . $$ The set $\La(\psi)$ is a
`$\limsup$' set; it consists of points in $\Omega$ which lie in
infinitely many of the balls $B(\ra,\p(\ma))$  centred at resonant
points. Clearly, even in this abstract setup  it is natural to refer
to the function $\p$ as the {\em approximating function}. It governs
the `rate' at which points in $\Omega$ must be approximated by
resonant sets in order to lie in $\La(\p)$.

Before continuing our discussion, we rewrite $\La(\p)$ in a
fashion which brings  its `$\limsup$' nature to the forefront. For
$n \in \N$, let $$
 \De(\p,n) :=  \!\!\!\!\!  \bigcup_{\alpha\in J \ : \ { k^{n-1} < \ma \leq k^n}}
\!\!\!\!\!\!\!\!\!\! B(\ra,\p(\ma)) \ \ \  \mbox{\rm where $k>1$
is fixed.} $$ By assumption the number of $ \alpha$ in $J$ with $
k^{n-1} < \ma \leq k^n $ is finite regardless of the value of $k$.
Thus, $\La(\psi)$ is precisely the set of points in $\Omega$ which
lie in infinitely many $\De(\p,n)$; that is $$ \La(\psi) =
\limsup_{n \to \infty} \De(\p,n) := \bigcap_{m=1}^{\infty}
\bigcup_{n=m}^{\infty} \De(\p,n) \ . $$

The main line of our investigation is motivated by the following
pair of fundamental  problems regarding the measure theoretic
structure of $\La(\p)$. In turn the fundamental problems are
motivated by the classical theory described in the  previous
section. It is reasonably straightforward to determine conditions
under which $m(\La(\p)) = 0$. In fact, this is implied by the
convergence part of the Borel--Cantelli lemma from probability
theory whenever
\begin{equation}\label{*1}
\sum_{n=1}^\infty m(\De(\p,n)) < \infty  \ .
\end{equation}
In view of this it is natural to consider:

\begin{problem}
Under what conditions is $m(\La(\p))$ strictly positive ?
\end{problem}

Under a `local ubiquity' hypothesis and a `$m$-volume' divergent sum
condition,  our first theorem provides a complete solution to this
problem; namely that $\La(\p)$ has full $m$--measure. This statement
can be viewed as the analogue of Khintchine's theorem.

Reiterating the above measure zero statement, if the approximating
function $\p$  decreases sufficiently quickly so that (\ref{*1})
is satisfied, the corresponding $\limsup $ set $\La(\p)$ is of
zero $m$--measure.  As with the classical setup of \S\ref{beg}, in
this case we cannot obtain any further information regarding the
`size' of $\La(\p)$ in terms of the ambient measure $m$  --- it is
always zero.  In short, we require a more delicate notion of
`size' than simply the given $m$-measure. In keeping with the
classical development, we investigate the `size' of $\La(\p)$ with
respect to  the Hausdorff measures $ \mathcal{ H}^{f} $ where $f$
is a dimension function.  Again, provided a certain `$f$-volume'
sum converges, it is reasonably simple to determine conditions
under which $\mathcal{ H}^{f}(\La(\p)) = 0$. Naturally, we
consider:

\begin{problem}
Under what conditions  is $\mathcal{ H}^{f}(\La(\p))$ strictly
positive ?
\end{problem}

 This problem turns out to be far more subtle  than
the previous one regarding $m$-measure. However, under a `local
ubiquity' hypothesis and an `$f$-volume' divergent sum condition,
together with mild conditions on the dimension function, our second
theorem shows that $\mathcal{ H}^{f}(\La(\p)) = \infty$. Thus,
$\mathcal{ H}^{f}(\La(\p))$ satisfies an  elegant `zero--infinity'
law whenever the convergence of the `$f$-volume' sum implies
$\mathcal{ H}^{f}(\La(\p)) = 0$ as is often the case. In particular,
this latter statement is true for the standard $s$-dimensional
Hausdorff measure $\mathcal{ H}^{s}$. Thus, in the language of
geometric measure theory the sets $\La(\p)$ are not $s$-sets.
Furthermore, from such zero--infinity laws it is easy to deduce the
Hausdorff dimension of $\La(\p)$.

In order to illustrate and clarify the above setup and our line of
investigation, we return to the  basic $\limsup$  set  of
\S\ref{beg}. The classical set $W(\p)$  of $\p$--well approximable
numbers in the theory of one dimensional Diophantine approximation
can clearly  be expressed in the form $\La(\p)$ with
\begin{eqnarray*}
& &\Omega:= [0,1] \  , \ \  J:= \{ (p,q) \in \N \times \N : 0\le p
\leq q \} \ , \ \ \alpha := (p,q) \in J \ , \
\\[1.5ex]
& &\ma := q \ , \ \ \ra:= p/q  \ \  \ {\rm and \ } \ \
\De(\ra,\p(\ma)) := B(p/q, \p(q)) \ .
\end{eqnarray*}
The metric $d$ is of course the standard Euclidean metric; $d(x,y)
:= |x-y| \, $.  Thus in this basic example, the resonant points
$\ra$ are simply rational points $p/q$. Furthermore,
$$\De(\p,n) := %[0.1] \, \cap
\bigcup_{k^{n-1} < q \leq k^n} \ \bigcup_{p=0}^{q} B(p/q, \p(q)) \
$$ and   $ W (\p) = \limsup \De(\p,n)$ as ${n \to \infty}$.

For this basic example, the solution to our first fundamental
problem is given by   Khintchine's theorem and the solution to the
second  by Jarn\'{\i}k's theorem. Together, these theorems
provide a complete measure theoretic description of $W(\p)$. In
the case of the general framework, analogues of these results
should be regarded as the ultimate  pair of results describing the
metric structure of the $\limsup$ sets $\La(\p)$. Alternatively,
they provide extremely satisfactory solutions  to the fundamental
problems. Analogues  of the convergence parts of the classical
results usually follow by adapting the `natural cover'
$$\{\De(\p,n) \, : \,  n=m, m+1, \cdots  \}   \hspace{16mm} (m \in
\N) \ $$ of $\La(\p)$. Our key  aim is to establish analogues of
the divergence parts of the classical results for $\La(\p)$.

%\section{Ubiquity and conditions on the general setup\label{congs}}

\section{Ubiquity \label{congs}}

In order to make any reasonable progress with the fundamental
problems we impose various conditions on the compact metric
measure space $(\Omega,d,m)$. Moreover, we require the notion of a
`local' ubiquitous system which will underpin our line of
investigation.

Throughout, a ball centred at a point $x$ and radius $r$ is
defined to be the set $\{ y \in \Omega : d(x,y) <  r \}$ or $\{ y
\in \Omega : d(x,y) \leq  r \}$ depending on whether it is open or
closed. In general,  we do not specify whether a certain ball is
open or close since it will be  irrelevant. Notice, that by
definition any ball is automatically a subset of $\Omega$. We
shall impose the following regularity condition on the measure of
balls.

\medskip

\noindent {\bf (M) \ } There exist positive constants $\delta$ and
$r_o$   such that for any $ x \in \Omega $ and $ r \leq r_o $,
\begin{equation*}
 a \, r ^{\delta}  \ \leq  \  m(B(x,r))  \ \leq  \   b \, r
^{\delta}    \ . \label{measure}
\end{equation*}

\medskip

\noindent The  constants $a$ and $b$ are independent of the ball
and without loss of generality we assume that $0<a<1<b$. Notice
that the above condition implies that $\dim \Omega = \delta $ and
furthermore that ${\cal H}^{\delta} (\Omega)$ is strictly positive
and finite. Indeed, $m$  is a comparable to $\delta$--dimensional
Hausdorff measure ${\cal H}^{\delta}$.

\subsection{The ubiquitous system\label{thesystems}}
The following `system' contain the key measure theoretic structure
necessary for our attack on the fundamental problems. Recall that
$\cR$ denotes the family of resonant sets $\ra$ and that the
function $\beta$ attaches a `weight' $\beta_{\alpha}$ to each
resonant set $\ra \in \cR$.

Let $\r : \R^+ \to \R^+ $ be a  function with $\r(r) \to 0 $ as $r
\to \infty $ and let

$$ \De(\r,n) := \bigcup_{\alpha\in J(n)} B(\ra,\r(k^n))
\ ,$$ where $k > 1 $ is a fixed real number and
$$
J(n) \ := \ \{ \alpha \in J : \beta_{\alpha} \leq k^n \} \ .
$$

\medskip

\noindent{\bf Definition (Local $\, m$--ubiquity)}~ Let $B=B(x,r)$
be an arbitrary ball with centre $x$ in $\Omega$ and radius $r \leq
r_0$. Suppose there exists a function $\r$  and  absolute constants
$\ka>0$ and $ k > 1$  such that
\begin{equation}\label{ub}
  m\left( B \cap \De(\r,n) \right) \ \ge \ \ka \ m(B) \qquad \mbox{for $n \geq n_o(B)$} .
\end{equation}
Then the pair $(\cR,\beta)$ is said to be a {\em local
$m$-ubiquitous system relative to $(\r,k)$}.

\medskip

Loosely speaking, the definition of local ubiquity says that the set
$\De(\r,n)$ locally approximates the underlying space $\Omega$ in
terms of the measure $m$. By `locally' we mean balls centred at
points in $\Omega$. The function $\r$, will be referred to as the
{\em ubiquitous function}. The actual values of the constants $\ka$
and $k$  in the above definition are irrelevant -- it is their
existence that is important.  In practice, the  local $m$--ubiquity
of a system can be established using standard arguments concerning
the distribution of the resonant sets in $\Omega$, from which the
function $\r$ arises naturally. To illustrate this, we return to the
classical  $\limsup$  set  of \S\ref{beg}.

The set $W(\p)$ of $\p$--well approximable numbers  has already been
shown to fit within our general $\limsup$ setup -- see
\S\ref{gensetup}. Now let $m$ be one--dimensional Lebesgue measure.
Clearly $m$  satisfies the measure condition (M) with $\delta =1 $.
With this in mind, we have the following statement concerning local
ubiquity within the classical setup.
%\vspace{1ex}
\begin{lemma}\label{reg_be}
There is a constant $k>1$ such that the pair $(\cR,\beta)$ is a
local $m$-ubiquitous system relative to $(\r,k)$ where  $\r : r
\mapsto {\rm constant} \times r^{-2}$.
\end{lemma}
\noindent{\em Proof.}~ \  Let $I=[a,b]\subset[0,1]$. By Dirichlet's
theorem
%\footnote{For any real number $x$ and integer $N \geq 1$,
%there exists coprime integers p,q  with $1\leq q \leq N $ such that
%$|x-p/q| < 1/(qN) \, . $}
, for any $x\in I$ there are coprime
integers $p,q$ with $1\le q\le k^n$ satisfying
$|x-p/q|<(qk^n)^{-1}$. Clearly, $aq-1\le p\le bq+1$. Thus, for a
fixed $q$ there are at most $m(I)q+3$ possible values of $p$.
Trivially, for $n$ large
 $$%\textstyle
m\left( I \cap  \!\! \bigcup_{q \le k^{n-1} }    \bigcup_{p} \
\textstyle{B\left(\frac{p}{q}, \frac{1}{q k^n} \right) } \right) \
\le \ 2 \!\! \sum_{q \, \leq k^{n-1} } \!\!\! \mbox{\large
$\frac{1}{qk^{n}} $ } (m(I)q+3) \le \mbox{\large $\frac{3}{k} $ }
 m(I)  \, .$$
 \noindent It follows that for  $k\geq 6$,
\begin{eqnarray*}
m\left( I \cap  \!\! \bigcup_{ q \le k^n } \bigcup_{p} \
\textstyle{B\left(\frac{p}{q}, \frac{k}{k^{2n}} \right) } \right)
 &  \geq &  m\left( I \cap  \!\! \bigcup_{k^{n-1} < q \le k^n }
\bigcup_{p} \ \textstyle{B\left(\frac{p}{q}, \frac{k}{k^{2n}}
\right) } \right) \\ & & \\  & \geq & \ m(I) - \mbox{\large
$\frac{3}{k} $ }m(I) \ \geq \ \mbox{\large $\frac{1}{2} $ } m(I) \
.
\end{eqnarray*}
%\end{proof}
 \QED

\medskip

It will be evident from our `ubiquity' theorems, that
Lemma~\ref{reg_be} is sufficient  for directly  establishing the
divergence part of both Khintchine's theorem and Jarn\'{\i}k's
zero--infinity law -- see \S\ref{basiceg}.

%%%%%%%%%%%%%%%%%%%%%%%%%%%%%%%%%%%%%%%%%%%%%%%%%%%%%%%%%%%%%%%%%%%%%%%%%%%%

%\subsection{A remark on related systems}

\vspace{3ex}

\noindent{\em A remark on related systems.} \ In the case that
$\Omega$ is a bounded subset of $\R^n$ and $m$ is $n$-dimensional
Lebesgue measure, the notion of ubiquity was originally formulated
by Dodson, Rynne \& Vickers \cite{DRV}  to obtain lower bounds for
the Hausdorff dimension of the sets $\Lambda(\p)$. Their notion of
ubiquity is closely related to our notion of a `local
$m$-ubiquitous' system and furthermore coincides with the `regular
systems' of Baker \& Schmidt \cite{BS}. Both these systems have
proved very useful in obtaining lower bounds for the Hausdorff
dimension of $\limsup$ sets. However, both \cite{BS} and
\cite{DRV} fail to shed any light on the problems considered in
this paper. For further details regarding regular systems and the
original formulation of ubiquitous systems see \cite{memoirs,BD}.
Recently and independently, in \cite{yann2} the notion of an
optimal regular system introduced in \cite{Ber2000} has been
re-formulated to obtain divergent type Hausdorff measures results
for subsets of $\R^n$. This re-formulated notion is essentially
equivalent to our notion of local $m$-ubiquity in which  $m$ is
$n$-dimensional Lebesgue measure and the ubiquity function is
comparable to $\rho: r \to r^{\mbox{\tiny $-1/n$ }} $.
Furthermore, even with these restrictions  our notion of local
$m$-ubiquity is not equivalent to that of an optimal regular
system since we make no assumption on the growth of $\# J(n)$.

%%%%%%%%%%%%%%%%%%%%%%%%%%%%%%%%%%%%%%%%%%%%%%%%%%%%%%%%%%%%%%%%%%%%%%%%%%

\subsection{The ubiquity statements} Recall, that an approximating
function $\psi$ is a real, positive decreasing function and that a
ubiquity function $\rho$ is a real, positive function such that
$\r(r) \to 0 $ as $r \to \infty $. Before stating our main results
we introduce one last notion. Given a real number $k
>1$, a function $h$ will be said to be {\bf ${\bf k}$-regular} if
there exists a strictly positive constant $\lambda < 1$ such that
for $n$ sufficiently large
\begin{equation} h(k^{n+1}) \leq \lambda \, h(k^n) \ .  \label{afmh}
\end{equation}
The constant $\lambda$ is independent of $n$ but may depend on
$k$. A consequence of local ubiquity are the following pair of
theorems. They constitute the main theorems appearing in
\cite{memoirs} tailored to the setup considered here.

\begin{thschmidt}
%\begin{theorem}
%\label{cor2} BDV1
Let $(\Omega,d)$ be a compact metric space equipped with a measure
$m$ satisfying condition\/ {\rm(M)}\/ such that any open  subset of
$\Omega$ is $m$--measurable . Suppose that $(\cR,\beta)$ is a local
$\,m$-ubiquitous system relative to $(\r,k)$ and that $\, \p$ is an
approximating function. Furthermore, suppose that either $\p$ or
$\r$ is $k$-regular and that
\begin{equation}  \label{dsm2}
  \sum_{n=1}^{\infty} \
 \left(\frac{\p(k^n)}{\r(k^n)}\right)^{\delta} \ = \ \infty  \  .
\end{equation}
 Then $$ m \left( \Lambda(\p) \right) = m(\Omega) \ .   $$
%\end{theorem}
\end{thschmidt}

\begin{thdet}
%\begin{theorem}\label{THM3}
 Let $(\Omega,d)$ be a compact metric space
equipped with a measure $m$ satisfying condition\/ {\rm(M)}\/.
Suppose that  $(\cR,\beta)$ is a locally $\,m$-ubiquitous system
relative
 to $(\rho,k)$ and that $\psi$ is an approximation
function. Let $f$ be a dimension function such that $r^{-\delta} \,
f(r) \to \infty$ as $r \to 0$ and $r^{-\delta} \, f(r)$  is
decreasing.  Let $g$ be the real, positive  function given by
\begin{equation}\label{my5}
g(r)  :=   f(\p(r))   \r(r)^{ -\delta} \ \mbox{ and let } \ \ G \,
:= \, \limsup_{n \to \infty} \, g(k^n).
\end{equation}
\begin{enumerate}
\item[(i)]
Suppose that $G = 0 $  and that $\r$ is $k$-regular. Then,
\begin{equation}\label{my6}
 {\cal H}^f \! \left( \Lambda(\p) \right)  =
 \infty \hspace{10mm} {\rm if } \hspace{10mm} \sum_{n=1}^{\infty}
g(k^n) =  \infty \; .
\end{equation}
\item[(ii)]
Suppose that $0< G \leq  \infty $. Then, ${\cal H}^f \! \left(
\Lambda(\p) \right) \ = \ \infty$.
\end{enumerate}
%\end{theorem}
\end{thdet}

Clearly, the assumption that the function $0< G \leq \infty $  in
part (ii) implies the divergent sum condition in part (i). The
case when the dimension function $f$ is $\delta$--dimensional
Hausdorff measure $\mathcal{ H}^{\delta}$ is excluded from the
statement of Theorem BDV2 by the condition that $r^{-\delta} \,
f(r) \to \infty$ as $r \to 0$. This is natural since otherwise
Theorem BDV1 implies that $m(\Lambda(\p))
> 0 $ which in turn implies that $\mathcal{ H}^{\delta} (\Lambda(\p))
$ is positive and finite.  In other words $\mathcal{ H}^{\delta}
(\Lambda(\p)) $ is never infinite.  However, given that the
measure $m$ is comparable to  $\mathcal{ H}^{\delta}$ -- a simple
consequence of condition (M) -- we are able to combine the above
statements and obtain  a {\bf single unifying theorem}.

\begin{theorem}\label{UNI}
 Let $(\Omega,d)$ be a compact metric space equipped with a measure
$m$ satisfying condition\/ {\rm(M)}\/ such that any open  subset
of $\Omega$ is $m$--measurable . Suppose that  $(\cR,\beta)$ is a
locally $\,m$-ubiquitous system relative
 to $(\rho,k)$ and that $\psi$ is an approximation
function. Let $f$ be a dimension function such that  $r^{-\delta}
\, f(r)$  is monotonic. Furthermore, suppose that  $\rho$ is
$k$-regular and that
\begin{equation}  \label{unif}
  \sum_{n=1}^{\infty} \
 \frac{f(\p(k^n))}{\r(k^n)^{\delta} }  \ = \ \infty  \  .
\end{equation}
Then,  $${\cal H}^f  \left( \Lambda(\p) \right) \ = \ {\cal H}^f
\left( \Omega \right) \ \ .
$$
\end{theorem}

The condition that $r^{-\delta} \, f(r)$  is monotonic is a
natural condition which is not particularly restrictive.  Note
that if the dimension function $f$ is such that  $r^{-\delta} \,
f(r) \to \infty$ as $r \to 0$ then $ {\cal H}^f \left( \Omega
\right)  = \infty $ and Theorem \ref{UNI} leads to the same
conclusion as Theorem BDV2. Here we make use of the following
fact: {\em if $ \, f$ and $g$ are two dimension functions such
that the ratio $f(r)/g(r) \to 0 $ as $ r \to 0 $, then $\mathcal{
H}^{f} (F) =0 $ whenever $\mathcal{ H}^{g} (F) < \infty $. }  On
the other hand, Theorem \ref{UNI} with $f(r) = r^{\delta}$ implies
that  ${\cal H}^{\delta} \left( \Lambda(\p) \right)  =  {\cal
H}^{\delta} \left( \Omega \right) $. This together with the fact
that the measure $m$ is comparable to ${\cal H}^{\delta}$ implies
that $m \left( \Lambda(\p) \right) = m(\Omega) $ -- the conclusion
of Theorem BDV1.

\subsection{The classical results\label{basiceg}}
For the classical set $W(\p)$ of $\p$--well approximable numbers,
 Lemma~\ref{reg_be} in \S\ref{congs}
establishes local $m$-ubiquity. Clearly, the ubiquity function
$\r$ satisfies (\ref{afmh}) (i.e. $\r$  is $u$-regular) and so
Theorem BDV1 establishes the divergent part of Khintchine's
Theorem. On the other hand, Theorem BDV2 establishes the divergent
part of Jarn\'{\i}k's Theorem. By making use of the `natural
cover' of $W(\p)$, the convergent parts of these classical results
are easily established.

In the above discussion we have opted to establish the classical
results of Khintchine and Jarn\'{\i}k separately. In the past
these results have always been thought of as separate entities
with Jarn\'{\i}k's  Theorem  being regarded as a refinement of
Khintchine's Theorem  -- but not containing Khintchine's Theorem.
However, it is easily seen that Lemma~\ref{reg_be} together with
Theorem \ref{UNI} leads to the following unification of the
fundamental classical results.

\begin{adeq}

\label{mainkj} Let $f$ be a dimension function such that   $r^{-1}
\, f(r) $ is monotonic. Let $\p$ be  a real, positive decreasing
function. Then $$ \hf\left(W(\p)\right)=\left\{\begin{array}{cl} 0
& {\rm \ if} \;\;\; \sum_{r=1}^{\infty}  \ \; r \,
f\left(\p(r)\right)
 <\infty \; ,\\[1ex]
\hf([0,1]) & {\rm \ if} \;\;\; \sum_{r=1}^{\infty} \  \; r \,
f\left(\p(r)\right) =  \infty           \; .
\end{array}\right.$$

\end{adeq}

\vspace{4ex}

\noindent{\em An important observation.}  It is worth standing
back a little and think about what we have actually used in
establishing the classical results -- namely local ubiquity.
Within the classical setup, local ubiquity is a simple measure
theoretic statement concerning the distribution of rational points
with respect to Lebesgue measure -- the natural measure on the
unit interval. From this we are able to obtain the divergent parts
of both Khintchine's Theorem (a Lebesgue measure statement) and
Jarn\'{\i}k's  Theorem (a Hausdorff measure statement).  In other
words, the Lebesgue measure statement of local ubiquity seems to
underpin the general Hausdorff measure theory of the $\limsup$ set
$W(\p)$.  That this is the case is by no means a coincidence --
see \cite{DS,DSslice}. In fact, in view of the Mass Transference
Principle introduced in \cite{DS} one actually has that
\begin{center} Khintchine's Theorem $ \hspace{4mm} \Longrightarrow
\hspace{4mm} $ Jarn\'{\i}k's Theorem.  \end{center} Thus, the
Lebesgue theory of $W(\p)$ underpins the general Hausdorff theory.
This at first glance is rather surprising in that the Hausdorff
theory had previously  been  thought to have been a subtle
refinement of the Lebesgue theory. However, given that the
Lebesgue statement of local ubiquity implies the general Hausdorff
theory we should not be too surprised.

\subsection{Where to go from here with ubiquity?  \label{desire}}
Let $\p$ be a real, positive decreasing function. For  $x$ in the
unit interval and $N \in \N$, let
$$
{\cal R}(x,N) \ := \ \# \{\; 1\leq q \leq N : | x - p/q|  < \p(q)
\; {\rm \ for \ some \ } p \in \Z \} \ .
$$
%Here $\|.\|$ denotes the distance to the nearest integer.
In view of Khintchine's Theorem, if $\sum q \p(q) $ diverges then
for almost all $x$ we have that ${\cal R}(x,N) \to \infty $ as $ N
\to \infty$. An obvious question now arises: can we saying
anything more precise about the behavior of the counting function
${\cal R}(x,N)$?  Within the classical theory of Diophantine
approximation we have the following remarkable quantitative
statement of Khintchine's Theorem.

\begin{thbs} Suppose that $2 \, q \, \psi(q) < 1$ and
that $\sum_{q=1}^{\infty} q \, \p(q) = \infty
$. Then, for almost all $x$
$$
{\cal R}(x,N) \ \sim \  \ 2 \,  \sum_{q=1}^{N} q \, \p(q) \ \ .
$$
\end{thbs}

\vspace{2ex}

\noindent Schmidt actually proves the above asymptotic statement
with an error term. Note that the condition $2 \, q \, \psi(q) <
1$ simply means that for any fixed $q$ there is at most one $p \in
\Z$ such that $ | x - p/q| < \p(q)$  -- this avoids counting
multiplicities.

In view of above discussion, in particular the work of Schmidt, a
gaping inadequacy with the ubiquity framework  is exposed. In
describing the $m$-measure theoretic structure of a $\limsup$, the
analogue of Khintchine's Theorem should be regarded as the first
step. The ultimate aim should be to obtain a quantitative version
of such a result; i.e. the analogue of Schmidt's Theorem. Thus we
ask the following question. Is there a natural `stronger' form of
local ubiquity which would enable us to obtain a quantitative form
of Theorem BDV1 analogous to Schmidt's Theorem? Obviously, it
would be highly desirable to establish such a form. Even a
ubiquity framework that would yield a comparable rather than
asymptotic analogue of Schmidt's Theorem would be desirable; i.e.
a framework which within the classical setup implies that $$ {\cal
R}(x,N) \ \asymp \   \sum_{q=1}^{N} q \, \p(q) \ .
$$

\section{Diophantine approximation and Kleinian Groups
\label{appkg}}
The classical results of Diophantine approximation,\hspace{2pt} in
particular those from the one dimensional theory, have natural
counterparts and extensions in the hyperbolic space setting. In
this setting, instead of approximating real numbers by rationals,
 one approximates limit points of a fixed   Kleinian group $G$ by
 points in the orbit (under the group) of a certain distinguished
 limit point $y$.
Beardon and Maskit  have shown that the geometry of the group is
reflected in the approximation properties of points in the limit
set.  The elements of $G$ are orientation preserving
 M\"obius transformations of the $(n+1)$--dimensional
 unit ball $B^{n+1}$. Let $\Lambda $ denote the limit set of $G$ and
 let $\delta$ denote the Hausdorff dimension of $\Lambda$.
 For any element $g$ in $G$ we shall use the notation
 $L_g := |g^\prime(0)|^{-1}$, where $|g^\prime(0)|$ is the (Euclidean)
 conformal dilation of $g$ at the origin.

Let $\psi $ be an approximating  function and let $$ W_{y}(\psi)
:=\{ \xi \in \Lambda: | \xi - g(y) | < \psi(L_g) \, \mbox{for i.m.
$g$ in $G$}\}  . $$

\noindent  This is the set of points in the limit set $\Lambda$
which are `close' to   infinitely  many  (`i.m.') images  of  a
`distinguished' point $y$. The `closeness' is of course governed
by the approximating function $\psi$. The limit point $y$ is taken
to be a parabolic fixed point if the group has parabolic elements
and a hyperbolic fixed point otherwise.

\vspace{2mm}

\noindent{\bf Geometrically finite groups with parabolics: \ } Let
us assume that the geometrically finite group has parabolic
elements so it is not convex co-compact. Thus our distinguished
limit point $y$ is a parabolic fixed point, say $p$. Associated
with $p$ is a geometrically motivated set $\mathcal{ T}_p$ of
coset representatives of $G_p \backslash G := \{g G_p : g \in G \}
$; so chosen that if $ g \in \mathcal{ T}_p$ then the orbit point
$g(0)$ of the origin lies within a bounded hyperbolic distance
from the top of the standard horoball $H_{g(p)}$. The latter, is
an $(n+1)$--dimensional Euclidean ball contained in $B^{n+1}$ such
that its boundary touches the unit ball $S^n$ at the point $g(p)$.
Let $R_g$ denote the Euclidean radius of  $H_{g(p)}$.  As a
consequence of the definition of $\mathcal{ T}_p$, it follows that
$$ \frac{1}{C \, L_g} \leq R_g \leq \frac{C}{L_g} \  $$ where
$C>1$ is an absolute constant. Also, it is worth mentioning that
the balls in the standard set of horoballs $\{H_{g(p)} \, : \,  g
\in \mathcal{ T}_p \} $ corresponding to the parabolic fixed point
$p$ are pairwise disjoint.   For further details and references
regarding the above notions and statements see any of the papers
\cite{jbgfg,melian,Strap}. With reference to our general
framework, let
 $\Omega:= \Lambda  \  , \  J:= \{ g: g \in \mathcal{ T}_p  \} \ , \ \alpha :=
g \in J \ , \ \ma := C \, R_g^{-1} \ {\rm and \ }
  \ra:=  g(p). $ Thus, the  family $\mathcal{ R}$ of resonant sets $\ra$ consists
  of orbit points $g(p)$
with $g \in \mathcal{ T}_p$. Furthermore, $B(\ra,\p(\ma)) :=
B(g(p),\p( C \, {R_g}^{-1} ) ) $ and
 $$\De(\p, n) :=   %\bigcup_{g \in J^*(n)}
 \bigcup_{\substack{g \in \mathcal{ T}_p \ :  \\  k^{n-1} < \,  C \, R_g^{-1} \leq
 k^n }}
 \!\!  B\left(g(p),\p( C \, R_g^{-1} ) \right) \ . $$
 %where $J^*(n) := \{ g \in \mathcal{ T}_p: k^{n-1} < \,  C \, R_g^{-1} \leq
 %k^n \} $.
 Here $k >1 $ is a constant.
Then  $$   W_{p}(\psi) \, \supset \, \Lambda(\p) \, := \,
\limsup_{n \to \infty} \De(\p, n) \ \ . $$

\noindent Now, let $m$  be  Patterson measure and  $ \delta = \dim
\Lambda$. Thus $m$ is a non-atomic, $\delta$--conformal
probability measure supported on $\Lambda$. We are assuming that
the group has parabolic elements, thus in general $m$ does not
satisfy condition (M) and so our ubiquity statements are not
applicable. However, if we restrict our attention to {\bf groups
of the first kind}  then $\Lambda = S^n$ and  $m$ is simply
$n$--dimensional Lebesgue measure on unit sphere $S^n$. Also note
that $\delta =n $ in  this case.  Thus for groups of the first
kind, $m$ clearly satisfies condition (M) and we have the
following statement concerning local ubiquity.

\medskip

\begin{proposition} Let $k \geq
k_o $ -- a positive group constant.  Then then pair $(\mathcal{
R}, \beta) $ is a local $m$--ubiquitous system relative to
$(\r,k)$ where $ \r : r \to \r(r) \, := \, {\rm constant} \ \times
\ r^{-1} \ . $ \label{kgps}
\end{proposition}

 The proposition follows from the following two facts which can
be found in  \cite{jbgfg,melian,Paddyrs}. They are valid in
general, but  for groups of the first kind they are particularly
easy to establish.

\smallskip

\noindent {\bf $\bullet$} {\em Local Horoball Counting Result: \ }
Let $B$ be an arbitrary Euclidean ball in $S^n$ centred at a limit
point. For $\lambda \in (0,1)$ and $r \in \R^+$ define $$
A_{\lambda}(B,R) \, := \, \{ g \in \mathcal{ T}_p: g(p) \in B {\rm \
and \ } \lambda R \leq R_g < R \}  \ . $$ There exists a positive
group constant $\lambda_o$ such that if $\lambda  \leq \lambda_o $
and $R < R_o(B)$, then

$$ k_1^{-1} \,  R^{-\delta} \,  m(B) \ \leq \ \# A_{\lambda}(B,R)
\ \leq \  k_1 \,  R^{-\delta} \,  m(B)  \ , $$ where $k_1$ is a
positive constant independent of $B$ and $R_o(B)$ is a
sufficiently small positive  constant  which does depend on $B$.

\noindent {\bf $\bullet$} {\em Disjointness Lemma: \ } For
distinct elements  $g,h \in \mathcal{ T}_p $  with $ \lambda <
R_g/R_h < \lambda^{-1}$, one has $ B(g(p),\lambda R_g ) \, \cap \,
B(h(p),\lambda R_h ) \ = \ \emptyset  \ . $

\vspace{2mm}

\noindent{\em Proof of Proposition \ref{kgps}.}~ \ To prove the
proposition, let $\rho(r) := C (k\, r)^{-1}$ where $k:= 1/\lambda
> 1/\lambda_o$ and $B$ be an arbitrary ball centred at a limit
point. Then for $n$ sufficiently large
\begin{eqnarray*}
m (\  B   \cap  \bigcup_{ g \in J(n) }
 B\left(g(p),\r( k^n ) \right)   \ )
& = & m (\  B   \cap  \bigcup_{\substack{ g \in \mathcal{ T}_p: \\
C \, R_g^{-1} \leq k^n} }
 B\left(g(p),\r( k^n ) \right)   \ )   \\ & & \\
 & \geq &
m  (   \bigcup_{\substack{g \in \mathcal{ T}_p: \, g(p) \in
\frac{1}{2} B \\
 k^{n-1} < C \, R_g^{-1} \leq k^n} }^{\circ}
 \!\!\!\!\!  B\left(g(p),\r( k^n ) \right)   \ )  \\ & &
\\
&\gg & m  (   \bigcup_{\substack{g \in \mathcal{ T}_p: \, g(p) \in
\frac{1}{2} B \\ k^{n-1} < \,
 C \, R_g^{-1} \leq k^n} }^{\circ}
 \!\!\!\!\!  B\left(g(p),\r( k^n ) \right)   \ )  \\ & &
\\
&\gg &  k^{-n \, \delta}  \ \#  A_{\frac{1}{k}}(\mbox{\small
$\frac{1}{2}$} B \, , C \, k^{-(n-1)} ) \\
&&\\
& \gg & m (\mbox{\small
$\frac{1}{2}$} B) \ \gg \ m(B) \  .
\end{eqnarray*}
\QED

\medskip

\noindent  Thus, in view of Proposition \ref{kgps} and the fact
that the measure $m$ is of type (M) and that $\r$ is $k$-regular,
Theorem \ref{UNI}  yields the divergent part of the following
statement. The convergent part is easy -- just use the `natural
cover' given by the $\limsup$ set $W_{p}(\psi)$  under
consideration. Also we make use of the following simple fact.
Suppose that $h: \R^+ \to \R^+ $ is a real, positive monotonic
function, $\alpha \in \R$ and $k>1$.  Then the divergence and
convergence properties of the sums $$ \sum_{n=1}^{\infty} k^{n \,
\alpha} \; h(k^n) \hspace{1cm} {\rm and }  \hspace{1cm}
\sum_{r=1}^{\infty} r^{ \alpha - 1 } \; h(r) \hspace{1cm} {\rm
coincide. }   $$

\medskip

\begin{theorem}\label{GF2}
Let $G$ be a geometrically finite Kleinian group of the first kind
with parabolic elements and $p$ be a parabolic fixed point.
 Let $f$ be a dimension function such that
 $r^{-n} \, f(r) $
is monotonic.  Let $\psi$  be a real,  positive decreasing
function. Then
$$ \hf\left( W_{p}(\psi)\right)= \left\{
\begin{array}{ll}
 0 & {\rm if}
 \;\;\; \sum_{r=1}^\infty \;  f\left(\p(r)\right) \;\;
 r^{n-1} <\infty\; ,\\[2ex]
\hf\left(S^n\right) & {\rm if} \;\;\; \sum_{r=1}^\infty \;
 f\left(\p(r)\right) \;\;  r^{n-1} =\infty \; .
\end{array}
\right. $$
\end{theorem}

\vspace{2ex}

In the above theorem, on taking $f(r) = r^n$ we obtain the
analogue of Khintchine's theorem with respect to the measure $m$
supported on the limit set; i.e. $n$--dimensional Lebesgue measure
on $ S^n$. The theorem is this case, under  a certain  regularity
condition on  $\psi$, has previously been established in
\cite{Paddyrs,Strap,SDS}. Indeed, in \cite{Strap} the analogue of
Khintchine's theorem with respect to  Patterson measure is
established without the condition that the group is of the first
kind. It is worth mentioning that the more general local ubiquity
framework of \cite{memoirs} also yields this statement even though
Patterson measure does not generically satisfy condition (M).
However, the condition (M) on the measure is essential even in
\cite{memoirs} for establishing general Hausdorff measure
`divergent' results and the full analogue of Theorem \ref{GF2}
without the `first kind' restriction is currently out of reach --
precise Hausdorff dimension statements are known \cite{jbgfg}.

%
% For example, in the case $\psi(r)
%= r^{-\tau}$ let us write $W_{p}(\tau) $ for $ W_{p}(\psi)$. Then
%$\dim W_{p}(\tau) = n/\tau $ $(\tau \geq 1)$. Clearly, Theorem
%\ref{GF2} implies this statement and shows that the
%$s$--dimensional Hausdorff measure of $W_{p}(\tau)$ at the
%critical exponent $s=n/\tau$ is infinite.

When interpreted on the upper half plane model $\Half^2$ of
hyperbolic space and applied to the modular group ${\rm
SL}(2,\Z)$, Theorem \ref{GF2} implies the classical result
associated with the  $\limsup$ set $W(\p)$ as stated in
\S\ref{basiceg}.

\medskip

 \noindent{\bf Convex co-compact groups:}~~~ These are geometrically
finite Kleinian groups without parabolic elements. Thus, the
distinguished limit point $y$ is a hyperbolic fixed point. For
convex co-compact groups, Patterson measure $m$ satisfies
condition (M) and the situation becomes much more satisfactory --
we don't not have to assume that the group is of the first kind.

 Let $L$ be the axis of the
conjugate pair of hyperbolic fixed points $y$ and $y'$, and let
$G_{yy'}$ denote the stabilizer of $y$ (or equivalently $y'$).
Then there is a geometrically motivated set $\mathcal{ T}_{yy'} $  of
coset representatives of $G_{yy'} \backslash G $; so chosen that
if $ g \in \mathcal{ T}_{yy'}$ then the orbit point $g(0)$ of the
origin lies within a bounded hyperbolic distance from the summit
$s_g$ of $g(L)$ -- the axis of the hyperbolic fixed pair $g(y)$
and $g(y')$.  The summit $s_g$ is simply the point on $g(L)$
`closest' to the origin. For $g \in \mathcal{ T}_{yy'}$, let
$H_{g(y)}$ be the horoball with base point at $g(y)$ and radius
$R_g := 1 - |s_g| $. Then the top of $H_{g(y)}$ lies within a
bounded hyperbolic distance of $g(0)$. Furthermore, as a
consequence of the definition of $\mathcal{ T}_{yy'}$, it follows that
$ C^{-1}  \leq R_g \, L_g \leq C $  where $C>1$ is an absolute
constant. We are now able to define the subset $\Lambda(\p)$ of
$W_{y}(\psi)$ in exactly the same way as in the parabolic  case
with $y$ replacing $p$ and $\mathcal{ T}_{yy'}$ replacing $\mathcal{
T}_{p}$.

Essentially the arguments given in \cite{melian}, can easily be
modified to obtain the analogue of the local  horoball counting
result stated above for the parabolic case.  We leave the details
to the reader. In turn, this enables  one to establish Proposition
\ref{kgps} for convex co-compact groups -- the statement remains
unchanged.  Since $m$ is of type (M) and $\rho$ is $k$--regular
for any $k>1$,  Theorem \ref{UNI} yields the divergent part of the
following statement. The convergent part is straightforward to
establish.

\medskip

\begin{theorem}
Let $G$ be a convex co-compact Kleinian group and $y$ be a
hyperbolic fixed point. Let $f$ be a dimension function such that
 $r^{-\delta}
\, f(r) $ is monotonic.  Let $\p $ be a real, positive decreasing
function. Then $$ \hf\left( W_{y}(\psi)\right)= \left\{
\begin{array}{ll}
  0 & {\rm if} \;\;\; \sum_{r=1}^\infty \;  f\left(\p(r)\right) \;\;
 r^{\delta-1} <\infty \; , \\[2ex]
 \hf\left(\Lambda\right) & {\rm if} \;\;\; \sum_{r=1}^\infty \;
 f\left(\p(r)\right) \;\;  r^{\delta-1} =\infty \; .
\end{array}\right.
$$ \label{CC2}
\end{theorem}

In the above theorem, on taking $f(r) = r^{\delta}$ we obtain the
convex co-compact analogue of Khintchine's theorem with respect to
the measure $m$ supported on the limit set; i.e. Patterson measure
on $\Lambda$. This Khintchine analogue, under  a certain regularity
condition on $\psi$, has been known for sometime -- see for example
\cite{gang4}. Regarding the general Hausdorff measure aspect of the
above theorem, previously only  dimension statements were known --
see \cite{slv3}. Theorem \ref{CC2} not only implies these dimension
statements but also gives the $s$--dimensional Hausdorff measure
${\cal H}^s$  of $W_{y}(\psi)$ at the critical exponent $s= \dim
W_{y}(\psi)$.

\subsection{Consequences of Theorem \ref{GF2}}

Throughout, $G$ is a  geometrically finite group of the first kind
with parabolic elements.  In this section we bring into play the
real strength of Theorem \ref{GF2}.  Let $\tau \geq 1 $ and
$\epsilon > 0 $ be arbitrary. Consider the approximating functions
$$ \p(r) \, :=  \, r^{-\tau} \, \left( \log \, r \right)^{
-  \frac{\tau}{n}  }  \hspace{12mm}  {\rm and  \ } \hspace{12mm}
\p_{\ep}(r) \, :=  \, r^{-\tau} \, \left( \log \, r \right)^{ -
\frac{\tau}{n} \left(1 + \ep \right) }  \; \; . $$ Let
$$ E_{p}(\tau) \ := \ W_{p}(\psi) \setminus W_{p}(\psi_\ep)  \ \ , $$
where $p$ is our `distinguished'  parabolic fixed point of $G$.
Thus, a limit point $\xi$ is in the set $ E(\tau) $ if
$$
| \xi - g(y) | < \psi(L_g) \hspace{6ex}  \mbox{for infinitely many
$g$ in $G$} \ ,
$$
and for any $\ep > 0$
$$
| \xi - g(y) | \geq \psi_\ep(L_g) \hspace{6ex}  \mbox{for all but
finitely many $g$ in $G$. }
$$
In other words, the approximation properties of $\xi$ by the orbit
of the parabolic fixed point  is `sandwiched' between the
approximating functions $\psi$ and $\psi_\ep$. Now consider the
dimension function
$$
f \, : \, r \to f(r) := r^{\frac{n}{\tau}} \  .
$$
Hence,  $\hf$ is simply $n/\tau$--dimensional Hausdorff measure
${\cal H}^{n/\tau}$. A straightforward application of Theorem
\ref{GF2} yields that
$$
{\cal H}^{\frac{n}{\tau}}(W_{p}(\p))  \, =  \,    {\cal
H}^{\frac{n}{\tau}}(S^n) \hspace{12mm} {\rm and  \ } \hspace{12mm}
{\cal H}^{\frac{n}{\tau}}(W_{p}(\psi_{\ep})) \, =  \,  0   \; \; .
$$ Now, $ {\cal H}^{n/\tau}(S^n)  > 0 $ (in fact it is equivalent to the
$n$--dimensional Lebesgue measure of the unit sphere $S^n$ when
$\tau = 1$ and is infinite if $\tau > 1$) and so we obtain the
following statement.

\begin{lemma} \label{exactlog} For $  \tau \geq 1 $, $ \dim
E_{p}(\tau) \, = \, n/\tau$ and furthermore
$$
{\cal H}^{\frac{n}{\tau}}(E_{p}(\tau))  \, =  \,    {\cal
H}^{\frac{n}{\tau}}(S^n)
$$
\end{lemma}

The main observation used in extracting Lemma \ref{exactlog} from
Theorem \ref{GF2} is the following: if we have two sets $A$ and
$B$ with $m(A) > 0 $ and $m(B) = 0 $ them $m(A \setminus B) = m(A)
> 0 $. This simply observation can be implemented  to obtain the
analogue of the lemma for general {\bf exact order} sets -- see
(\cite{BS,exactus3,Bugbook,Gut} for a discussion of this notion
within the classical framework of Diophantine approximation.
Briefly, given two approximating functions $\varphi$ and $\psi$
with $\varphi$ in some sense `smaller' than $\psi$, consider the
set $E_p(\psi,\varphi) := W_p(\psi) \setminus W_p(\varphi)$. Thus
the approximation properties of limit points $\xi$ in
$E_p(\psi,\varphi)$ are `sandwiched' between the functions
$\varphi$ and $\psi$. Under suitable conditions on the `smallness'
of $\varphi$ compared to $\psi$ it is possible to obtain the
analogue of Lemma \ref{exactlog} for the set $E_p(\psi,\varphi)$
-- see \cite{exactus3} for the classical statements. In view of
the above observation,  the key is to construct an appropriate
dimension function $f$ for which $\mathcal{ H}^f(W_p(\psi))=
\hf(S^n)  $ and $\mathcal{ H}^f(W_p(\varphi))=0$.

\vspace{2ex}

In the case that $\tau =1$,  Lemma \ref{exactlog} has a well known
dynamical interpretation in terms of the `rate' of excursions by
geodesics into a cuspidal end of the associated hyperbolic
manifold  $ {\cal M} = B^{n+1}/G $; namely Sullivan's logarithm
law for geodesics \cite{SDS}. We are now in the position to
naturally place this law within the general Hausdorff measure
setting. First some notation.  Let $P$ denote a complete set of
parabolic fixed points inequivalent under $G$. Clearly the orbit
$G(P)$ of points in $P$ under $G$ is the complete set of parabolic
fixed points of $G$.  Since $G$ is geometrically finite of the
first kind with parabolic elements, the associated hyperbolic
manifold $ {\cal M} $ consists of a compact part $ X_o $ with a
finite number of attachments:
$$
{\cal M} =  X_o \; \; \cup  \; \;  \bigcup_{p \in P} Y_p  %\; \;
%\cup \; \;
% \bigcup_{i=1}^{k} Z_i \;
$$
where each $p$ in $P$ determines an exponentially `thinning' end $
Y_{p} $ -- usually referred to as a cuspidal end -- attached to
$X_o$.  %and the free faces of a convex fundamental polyhedron for
%$G$ correspond to a finite collection of exponentially `exploding'
%ends $Z_i$ (usually referred to as a funnels) attached to $X_o \,$
%-- see \cite{Nic,SDS}. In the situation that the group is of the
%first kind, as is the case here, there are no funnels.

We shall write $0$ for the projection of the origin in $B^{n+1}$
to the quotient space $\cal M$. Let $S^n$ be the unit sphere of
the tangent space to $\cal M$ at $0$, and for every vector $v$ in
$S^n$ let $\gamma_v$ be the geodesic emanating from $0$ in the
direction $v$. Furthermore, for $t$ in $\R^+$, let $ \gamma_v (t)
$ denote the point achieved after travelling time $t$ along
$\gamma_v$. Now fix a $p\in P$. We define a function
\begin{eqnarray*}
 \pen_{p}:{\cal M} & \to & \R^+\\
  x &\mapsto &\left\{\begin{array}{ll}
  0 & x\notin Y_{p}\\
  \dist(x,0) & x\in Y_{p},
  \end{array}\right.
\end{eqnarray*}
where dist is the induced metric on $\cal M$. This is the
penetration of $x$ into the cuspidal end $Y_{p}$.  A relatively
simple argument (see \cite{MP,SDS,slv3}) shows that the excursion
pattern of a random geodesic into a cuspidal end $Y_p$ is
equivalent to the approximation of a random limit point of $G$ by
the base points of standard horoballs in $\{H_{g(p)} \, : \,  g
\in \mathcal{ T}_p \} $. In particular, for any $ \alpha $ in
$[0,1] $, consider the set ${\cal S}_p(\alpha) $ of directions $v$
in $S^n$ such that
$$
\limsup_{t \to \infty} \frac{ \pen_{p}(\gamma_v ( t))  \, - \,
\alpha t }{ \log t } \  =  \  \frac{1}{n}  \ \ .
$$
Then the problem of determining the measure theoretic structure of
${\cal S}_p(\alpha) $ is equivalent to  determining the measure
theoretic structure of $ E_{p}(\tau)$ with  $\tau = 1/(1-\alpha)$.
In view of this, the following result can be regarded as a
dynamical interpretation of Lemma \ref{exactlog} in terms of the
geodesic excursions into the cuspidal ends of ${\cal M}$.

\begin{theorem}[A general logarithm law for geodesics]
\label{gsul} Let $G$ be a geometrically finite group of the first
kind with parabolic elements. For $\alpha \in [0,1)$, we have that
$$
{\cal H}^{n(1-\alpha)}( {\cal S}_p(\alpha) )  \; = \; {\cal
H}^{n(1-\alpha)}(S^n) \;.
$$
\end{theorem}

 In the case $\alpha = 0$, so that ${\cal H}^{n(1-\alpha)}$ is equivalent to
$n$-dimensional Lebesgue measure, the theorem reduces to
Sullivan's famous logarithm law for geodesics. The theorem simple
says that Sullivan's logarithm law survives for $ \alpha > 0 $ if
we appropriately `rescale' $n$-dimensional Lebesgue measure.

\medskip

\noindent{\bf Remark.} In this section we have chosen to
demonstrate the power of Theorem \ref{GF2}. We could just as
easily have  picked on Theorem \ref{CC2} and established analogues
statements to Lemma \ref{exactlog}  and Theorem \ref{gsul} for
convex co-compact groups. The latter would be a statement along
the lines suggested by the dynamical interpretation of the
Diophantine approximation results in \cite{gang4}.

\medskip

We end our discussion by studying limit points which are
`extremely' well approximable by the orbit of a parabolic fixed
point. In view of the above discussion, they correspond to
geodesics which exhibit an `extremely' rapid excursion pattern
into a cuspidal end of ${\cal M}$. For $\omega
> 0 $, let us say that a limit point $ \xi$ is {\bf
$\omega$--Liouville} if
$$
| \xi - g(y) | < \exp(-L_g^{\omega}) \hspace{3ex} \mbox{for
infinitely many $g$ in $G$}  \ .
$$
Let $L_p(\omega)$ denote the set of $\omega$--Liouville limit
points. Note that if $\xi \in L_p(\omega)$, then for any real
number $\tau$ we have that  $| \xi - g(y) | < L_g^{-\tau}$ for
infinitely many  $g$ in $G$ -- hence the reference to Liouville
since in the classical framework, a real number $x$ is said to be
Liouville if $|x - p/q | < q^{-\tau}$ for infinitely many
rationals $p/q$, irrespective of the value of $\tau$. It is easy
to see that for any $s > 0$
$$
\sum_{r=1}^{\infty} r^{n-1} \left( \exp(-r^{\omega})\right)^s  \ \
< \ \infty  \ , $$ regardless of $\omega$ and so the sets
$L_p(\omega)$ are of zero dimension.
%Thus the usual notion of
%Hausdorff dimension is unable to distinguish between the sets
However, given $\ep \geq 0$, let $f_{\ep}$  be the dimension
function given by
$$
%f(r)  :=  \left(\log \frac{1}{r}\right)^{\frac{n}{\omega}} \times
%\left(\log \log \frac{1}{r}\right)^{-1}  \hspace{2ex} {\rm and }
%\hspace{2ex}
f_{\ep}(r)  :=  \left(\log
\frac{1}{r}\right)^{\frac{n}{\omega}} \times \left(\log \log
\frac{1}{r}\right)^{-(1+\ep)} \ .
$$
On applying Theorem \ref{GF2}, we obtain the following statement.

\begin{lemma}
\label{love} Let $G$ be a geometrically finite group of the first
kind with parabolic elements. For $ \omega > 0 $,
$$ {\cal H}^{f_{\ep}} (L_p(\omega)) = \left\{
\begin{array}{ll}
  0 & {\rm if} \;\;\; \ep > 0 \; , \\[2ex]
 \infty & {\rm if} \;\;\; \ep = 0 \; .
\end{array}\right.
$$
\end{lemma}

\vspace{3ex}

In terms of the dimension theory,  when we are confronted with
sets of dimension zero it is natural to change the usual
`$r^s$-scale' in the definition of Hausdorff dimension to a
logarithmic scale. For $s
> 0$, let $f_s$ be the dimension function given by $f_s(r):=
(-\log r)^s$. The {\bf logarithmic Hausdorff dimension} of a set
$F$ is defined by $ \dim_{\substack{ \\   \log}}  F :=  \inf
\left\{ s : \mathcal{ H}^{f_s} (F) =0 \right\} = \sup \left\{ s :
\mathcal{ H}^{f_s} (F) = \infty \right\} $. It is easily verified
that if $ \dim F > 0 $ then $ \dim_{\substack{ \\   \log}}  F =
\infty$ -- precisely as one should expect.  The following
statement %concerning the logarithmic dimension of $L_p(\omega)$
is
a simple consequence of Lemma \ref{love}.

\begin{corollary}
\label{corlove} Let $G$ be a geometrically finite group of the
first kind with parabolic elements. For $ \omega > 0 $,
$$
\dim_{\substack{ \\   \log}}  L_p(\omega) = \frac{n}{\omega}  \ .
$$
Furthermore, ${\cal H}^{f_{s}} (L_p(\omega)) = \infty $ at the
critical exponent $s= n/\omega$.
\end{corollary}

\medskip

\noindent{\bf Remark.} Theorem \ref{CC2} yields the analogues
statements to Lemma \ref{love}  and Corollary \ref{corlove} for
convex co-compact groups. Apart from replacing $n$ by $\delta :=
\dim \Lambda $, the statements are identical to those above.

\vspace{6mm}

\noindent {\bf Acknowledgements.} SV would like to thank Francoise
Dalbo and Cornelia Drutu for organising the conference on `Dynamical
systems and Diophantine Approximation', held at the Institut Henri
Poincaré, Paris, 7-9 June 2004 and for giving me the opportunity to
give a series of lectures -- it was a most enjoyable experience.
Also a special thanks to Francoise for putting up with me in Rennes
after the conference -- it could not have been easy given that she
was just about to give birth to twins. I still savour the taste of
those excellent crepes! Finally, I would like to thank my pair --
Ayesha and Iona -- for nearly five wonderful years of tear-producing
laughter.  As they often remind me, I am still the Scarecrow in
search of his brain -- thanks a bunch girls!

\vspace{5mm}

\noindent Victor V. Beresnevich: Department of Mathematics,
University of York,

\vspace{-2mm}

\noindent\phantom{Victor V. Beresnevich: }Heslington, York, YO10
5DD, England.

%\vspace{0mm}

\noindent\phantom{Victor V. Beresnevich: }e-mail: vb8@york.ac.uk

\vspace{5mm}

\noindent Sanju L. Velani: Department of Mathematics, University
of York,

\vspace{-2mm}

 ~ \hspace{19mm}  Heslington, York, YO10 5DD, England.

%\vspace{0mm}

 ~ \hspace{19mm} e-mail: slv3@york.ac.uk

\end{document}